\documentclass{amsart}
\usepackage{amsfonts}

\newtheorem{theorem}{Theorem}
\theoremstyle{plain}

\newtheorem{corollary}{Corollary}

\newtheorem{lemma}{Lemma}

\newtheorem{proposition}{Proposition}
\newtheorem{remark}{Remark}

\numberwithin{equation}{section}
\input{tcilatex}

\begin{document}
\title[Quadratic Reverses of the Triangle Inequality]{Quadratic Reverses of
the Triangle Inequality in Inner Product Spaces}
\author{Sever S. Dragomir}
\address{School of Computer Science and Mathematics\\
Victoria University of Technology\\
PO Box 14428, MCMC 8001\\
Victoria, Australia.}
\email{sever@csm.vu.edu.au}
\urladdr{http://rgmia.vu.edu.au/SSDragomirWeb.html}
\date{April 16, 2004.}

\begin{abstract}
Some sharp quadratic reverses for the generalised triangle inequality in
inner product spaces and applications are given.
\end{abstract}

\maketitle

\section{Introduction}

In 1966, J.B. Diaz and F.T. Metcalf \cite{DM} proved the following reverse
of the triangle inequality in the general settings of inner product spaces:

\begin{theorem}
\label{ta}Let $a$ be a unit vector in the inner product space $\left(
H;\left\langle \cdot ,\cdot \right\rangle \right) $ over the real or complex
number field $\mathbb{K}$. Suppose that the vectors $x_{i}\in H\backslash
\left\{ 0\right\} ,$ $i\in \left\{ 1,\dots ,n\right\} $ satisfy%
\begin{equation}
0\leq r\leq \frac{\func{Re}\left\langle x_{i},a\right\rangle }{\left\Vert
x_{i}\right\Vert },\ \ \ \ \ i\in \left\{ 1,\dots ,n\right\} .  \label{1.1}
\end{equation}%
Then%
\begin{equation}
r\sum_{i=1}^{n}\left\Vert x_{i}\right\Vert \leq \left\Vert
\sum_{i=1}^{n}x_{i}\right\Vert ,  \label{1.2}
\end{equation}%
where equality holds if and only if%
\begin{equation}
\sum_{i=1}^{n}x_{i}=r\left( \sum_{i=1}^{n}\left\Vert x_{i}\right\Vert
\right) a.  \label{1.3}
\end{equation}
\end{theorem}

For some similar results valid for semi-inner products in normed spaces, see 
\cite{K} and \cite{M}.

In the same spirit, but providing a somewhat simpler sufficient condition
with a clear geometrical meaning, we note the following result obtained by
the author in \cite{SSD1}:

\begin{theorem}
\label{tb}Let $a$ be as above and $\rho \in \left( 0,1\right) .$ If $%
x_{i}\in H,$ $i\in \left\{ 1,\dots ,n\right\} $ are such that%
\begin{equation}
\left\Vert x_{i}-a\right\Vert \leq \rho \text{ \ for each \ }i\in \left\{
1,\dots ,n\right\} ,  \label{1.4}
\end{equation}%
then we have the inequality%
\begin{equation}
\sqrt{1-\rho ^{2}}\sum_{i=1}^{n}\left\Vert x_{i}\right\Vert \leq \left\Vert
\sum_{i=1}^{n}x_{i}\right\Vert ,  \label{1.5}
\end{equation}%
with equality if and only if%
\begin{equation}
\sum_{i=1}^{n}x_{i}=\sqrt{1-\rho ^{2}}\left( \sum_{i=1}^{n}\left\Vert
x_{i}\right\Vert \right) a.  \label{1.6}
\end{equation}
\end{theorem}

In a complementary direction providing reverses of the triangle inequality
in its additive form, i.e., upper bounds for the nonnegative difference%
\begin{equation*}
\sum_{i=1}^{n}\left\Vert x_{i}\right\Vert -\left\Vert
\sum_{i=1}^{n}x_{i}\right\Vert ,
\end{equation*}%
we note the following recent result obtained in \cite{SSD1}:

\begin{theorem}
\label{tc}Let $a$ be as above and $x_{i}\in H,$ $k_{i}\geq 0,$ $i\in \left\{
1,\dots ,n\right\} $ such that 
\begin{equation}
\left\Vert x_{i}\right\Vert -\func{Re}\left\langle a,x_{i}\right\rangle \leq
k_{i}\text{ \ for each \ }i\in \left\{ 1,\dots ,n\right\} ,  \label{1.7}
\end{equation}%
then we have the inequality%
\begin{equation}
0\leq \sum_{i=1}^{n}\left\Vert x_{i}\right\Vert -\left\Vert
\sum_{i=1}^{n}x_{i}\right\Vert \leq \sum_{i=1}^{n}k_{i}.  \label{1.8}
\end{equation}%
The equality holds in (\ref{1.8}) if and only if%
\begin{equation}
\sum_{i=1}^{n}\left\Vert x_{i}\right\Vert \geq \sum_{i=1}^{n}k_{i}
\label{1.9}
\end{equation}%
and 
\begin{equation}
\sum_{i=1}^{n}x_{i}=\left( \sum_{i=1}^{n}\left\Vert x_{i}\right\Vert
-\sum_{i=1}^{n}k_{i}\right) a.  \label{1.10}
\end{equation}
\end{theorem}

Another similar result but with a simpler condition, is the following one 
\cite{SSD1}.

\begin{theorem}
\label{td}Let $a$ and $x_{i},$ $i\in \left\{ 1,\dots ,n\right\} $ be as
above. If $r_{i}>0,$ $i\in \left\{ 1,\dots ,n\right\} $ are such that%
\begin{equation}
\left\Vert x_{i}-a\right\Vert \leq r_{i}\text{ \ for each \ }i\in \left\{
1,\dots ,n\right\} ,  \label{1.11}
\end{equation}%
then we have the inequality%
\begin{equation}
0\leq \sum_{i=1}^{n}\left\Vert x_{i}\right\Vert -\left\Vert
\sum_{i=1}^{n}x_{i}\right\Vert \leq \frac{1}{2}\sum_{i=1}^{n}r_{i}^{2}.
\label{1.12}
\end{equation}%
The equality holds in (\ref{1.12}) if and only if%
\begin{equation}
\sum_{i=1}^{n}\left\Vert x_{i}\right\Vert \geq \frac{1}{2}%
\sum_{i=1}^{n}r_{i}^{2}  \label{1.13}
\end{equation}%
and%
\begin{equation}
\sum_{i=1}^{n}x_{i}=\left( \sum_{i=1}^{n}\left\Vert x_{i}\right\Vert -\frac{1%
}{2}\sum_{i=1}^{n}r_{i}^{2}\right) a.  \label{1.14}
\end{equation}
\end{theorem}

For other inequalities related to the triangle inequality, see Chapter XVII
of the book \cite{MPF}.

The main aim of the present paper is to point out some quadratic reverses
for the generalised triangle inequality, namely, sharp upper bounds for the
nonnegative differences%
\begin{equation*}
\left( \sum_{i=1}^{n}\left\Vert x_{i}\right\Vert \right) ^{2}-\left\Vert
\sum_{i=1}^{n}x_{i}\right\Vert ^{2},
\end{equation*}%
under various assumptions for the vectors $x_{i}\in H,$ $i\in \left\{
1,\dots ,n\right\} $ involved. Some related results are established.
Applications for vector-valued integrals in Hilbert spaces are also given.

\section{Quadratic Reverses of the Triangle Inequality}

The following lemma holds:

\begin{lemma}
\label{l2.1}Let $\left( H;\left\langle \cdot ,\cdot \right\rangle \right) $
be an inner product space over the real or complex number field $\mathbb{K}$%
, $x_{i}\in H,$ $i\in \left\{ 1,\dots ,n\right\} $ and $k_{ij}>0$ for $1\leq
i<j\leq n$ such that%
\begin{equation}
0\leq \left\Vert x_{i}\right\Vert \left\Vert x_{j}\right\Vert -\func{Re}%
\left\langle x_{i},x_{j}\right\rangle \leq k_{ij}  \label{2.1}
\end{equation}%
for $1\leq i<j\leq n.$ Then we have the following quadratic reverse of the
triangle inequality%
\begin{equation}
\left( \sum_{i=1}^{n}\left\Vert x_{i}\right\Vert \right) ^{2}\leq \left\Vert
\sum_{i=1}^{n}x_{i}\right\Vert ^{2}+2\sum_{1\leq i<j\leq n}k_{ij}.
\label{2.2}
\end{equation}%
The case of equality holds in (\ref{2.2}) if and only if it holds in (\ref%
{2.1}) for each $i,j$ with $1\leq i<j\leq n.$
\end{lemma}

\begin{proof}
We observe that the following identity holds:%
\begin{align}
& \left( \sum_{i=1}^{n}\left\Vert x_{i}\right\Vert \right) ^{2}-\left\Vert
\sum_{i=1}^{n}x_{i}\right\Vert ^{2}  \label{2.3} \\
& =\sum_{i,j=1}^{n}\left\Vert x_{i}\right\Vert \left\Vert x_{j}\right\Vert
-\left\langle \sum_{i=1}^{n}x_{i},\sum_{j=1}^{n}x_{j}\right\rangle  \notag \\
& =\sum_{i,j=1}^{n}\left\Vert x_{i}\right\Vert \left\Vert x_{j}\right\Vert
-\sum_{i,j=1}^{n}\func{Re}\left\langle x_{i},x_{j}\right\rangle  \notag \\
& =\sum_{i,j=1}^{n}\left[ \left\Vert x_{i}\right\Vert \left\Vert
x_{j}\right\Vert -\func{Re}\left\langle x_{i},x_{j}\right\rangle \right] 
\notag \\
& =\sum_{1\leq i<j\leq n}\left[ \left\Vert x_{i}\right\Vert \left\Vert
x_{j}\right\Vert -\func{Re}\left\langle x_{i},x_{j}\right\rangle \right]
+\sum_{1\leq j<i\leq n}\left[ \left\Vert x_{i}\right\Vert \left\Vert
x_{j}\right\Vert -\func{Re}\left\langle x_{i},x_{j}\right\rangle \right] 
\notag \\
& =2\sum_{1\leq i<j\leq n}\left[ \left\Vert x_{i}\right\Vert \left\Vert
x_{j}\right\Vert -\func{Re}\left\langle x_{i},x_{j}\right\rangle \right] . 
\notag
\end{align}%
Using the condition (\ref{2.1}), we deduce that%
\begin{equation*}
\sum_{1\leq i<j\leq n}\left[ \left\Vert x_{i}\right\Vert \left\Vert
x_{j}\right\Vert -\func{Re}\left\langle x_{i},x_{j}\right\rangle \right]
\leq \sum_{1\leq i<j\leq n}k_{ij},
\end{equation*}%
and by (\ref{2.3}), we deduce the desired inequality (\ref{2.2}).

The case of equality is obvious by the identity (\ref{2.3}) and we omit the
details.
\end{proof}

\begin{remark}
\label{r2.2}From (\ref{2.2}) one may deduce the coarser inequality that
might be useful in some applications:%
\begin{equation*}
0\leq \sum_{i=1}^{n}\left\Vert x_{i}\right\Vert -\left\Vert
\sum_{i=1}^{n}x_{i}\right\Vert \leq \sqrt{2}\left( \sum_{1\leq i<j\leq
n}k_{ij}\right) ^{\frac{1}{2}}\quad \left( \leq \sqrt{2}\sum_{1\leq i<j\leq
n}\sqrt{k_{ij}}\right) .
\end{equation*}
\end{remark}

\begin{remark}
\label{r2.3}If the condition (\ref{2.1}) is replaced with the following
refinement of Schwarz's inequality:%
\begin{equation}
\left( 0\leq \right) \delta _{ij}\leq \left\Vert x_{i}\right\Vert \left\Vert
x_{j}\right\Vert -\func{Re}\left\langle x_{i},x_{j}\right\rangle \text{ for }%
1\leq i<j\leq n,  \label{2.4}
\end{equation}%
then the following refinement of the quadratic generalised triangle
inequality is valid:%
\begin{equation}
\left( \sum_{i=1}^{n}\left\Vert x_{i}\right\Vert \right) ^{2}\geq \left\Vert
\sum_{i=1}^{n}x_{i}\right\Vert ^{2}+2\sum_{1\leq i<j\leq n}\delta _{ij}\quad
\left( \geq \left\Vert \sum_{i=1}^{n}x_{i}\right\Vert ^{2}\right) .
\label{2.5}
\end{equation}%
The equality holds in the first part of (\ref{2.5}) iff the case of equality
holds in (\ref{2.4}) for each $1\leq i<j\leq n.$
\end{remark}

The following result holds.

\begin{proposition}
\label{p2.4}Let $\left( H;\left\langle \cdot ,\cdot \right\rangle \right) $
be as above, $x_{i}\in H,$ $i\in \left\{ 1,\dots ,n\right\} $ and $r>0$ such
that%
\begin{equation}
\left\Vert x_{i}-x_{j}\right\Vert \leq r  \label{2.6}
\end{equation}%
for $1\leq i<j\leq n.$ Then%
\begin{equation}
\left( \sum_{i=1}^{n}\left\Vert x_{i}\right\Vert \right) ^{2}\leq \left\Vert
\sum_{i=1}^{n}x_{i}\right\Vert ^{2}+\frac{n\left( n-1\right) }{2}r^{2}.
\label{2.7}
\end{equation}%
The case of equality holds in (\ref{2.7}) if and only if 
\begin{equation}
\left\Vert x_{i}\right\Vert \left\Vert x_{j}\right\Vert -\func{Re}%
\left\langle x_{i},x_{j}\right\rangle =\frac{1}{2}r^{2}  \label{2.8}
\end{equation}%
for each $i,j$ with $1\leq i<j\leq n.$
\end{proposition}

\begin{proof}
The inequality (\ref{2.6}) is obviously equivalent to%
\begin{equation*}
\left\Vert x_{i}\right\Vert ^{2}+\left\Vert x_{j}\right\Vert ^{2}\leq 2\func{%
Re}\left\langle x_{i},x_{j}\right\rangle +r^{2}
\end{equation*}%
for $1\leq i<j\leq n.$ Since%
\begin{equation*}
2\left\Vert x_{i}\right\Vert \left\Vert x_{j}\right\Vert \leq \left\Vert
x_{i}\right\Vert ^{2}+\left\Vert x_{j}\right\Vert ^{2},\ \ 1\leq i<j\leq n;
\end{equation*}%
hence%
\begin{equation}
\left\Vert x_{i}\right\Vert \left\Vert x_{j}\right\Vert -\func{Re}%
\left\langle x_{i},x_{j}\right\rangle \leq \frac{1}{2}r^{2}  \label{2.9}
\end{equation}%
for any $i,j$ with $1\leq i<j\leq n.$

Applying Lemma \ref{l2.1} for $k_{ij}:=\frac{1}{2}r^{2}$ and taking into
account that 
\begin{equation*}
\sum_{1\leq i<j\leq n}k_{ij}=\frac{n\left( n-1\right) }{4}r^{2},
\end{equation*}
we deduce the desired inequality (\ref{2.7}). The case of equality is also
obvious by the above lemma and we omit the details.
\end{proof}

In the same spirit, and if some information about the forward difference $%
\Delta x_{k}:=x_{k+1}-x_{k}$ $\left( 1\leq k\leq n-1\right) $ are available,
then the following simple quadratic reverse of the generalised triangle
inequality may be stated.

\begin{corollary}
\label{c2.5}Let $\left( H;\left\langle \cdot ,\cdot \right\rangle \right) $
be an inner product space and $x_{i}\in H,$ $i\in \left\{ 1,\dots ,n\right\}
.$ Then we have the inequality%
\begin{equation}
\left( \sum_{i=1}^{n}\left\Vert x_{i}\right\Vert \right) ^{2}\leq \left\Vert
\sum_{i=1}^{n}x_{i}\right\Vert ^{2}+\frac{n\left( n-1\right) }{2}%
\sum_{k=1}^{n-1}\left\Vert \Delta x_{k}\right\Vert .  \label{2.10}
\end{equation}%
The constant $\frac{1}{2}$ is best possible in the sense that it cannot be
replaced in general by a smaller quantity.
\end{corollary}

\begin{proof}
Let $1\leq i<j\leq n.$ Then, obviously, 
\begin{equation*}
\left\Vert x_{j}-x_{i}\right\Vert =\left\Vert \sum_{k=i}^{j-1}\Delta
x_{k}\right\Vert \leq \sum_{k=i}^{j-1}\left\Vert \Delta x_{k}\right\Vert
\leq \sum_{k=1}^{n-1}\left\Vert \Delta x_{k}\right\Vert .
\end{equation*}%
Applying Proposition \ref{p2.4} for $r:=\sum_{k=1}^{n-1}\left\Vert \Delta
x_{k}\right\Vert ,$ we deduce the desired result (\ref{2.10}).

To prove the sharpness of the constant $\frac{1}{2},$ assume that the
inequality (\ref{2.10}) holds with a constant $c>0,$ i.e., 
\begin{equation}
\left( \sum_{i=1}^{n}\left\Vert x_{i}\right\Vert \right) ^{2}\leq \left\Vert
\sum_{i=1}^{n}x_{i}\right\Vert ^{2}+cn\left( n-1\right)
\sum_{k=1}^{n-1}\left\Vert \Delta x_{k}\right\Vert  \label{2.11}
\end{equation}%
for $n\geq 2,$ $x_{i}\in H,$ $i\in \left\{ 1,\dots ,n\right\} .$

If we choose in (\ref{2.11}), $n=2,$ $x_{1}=-\frac{1}{2}e,$ $x_{2}=\frac{1}{2%
}e,$ $e\in H,$ $\left\Vert e\right\Vert =1,$ then we get $1\leq 2c,$ giving $%
c\geq \frac{1}{2}.$
\end{proof}

The following result providing a reverse of the quadratic generalised
triangle inequality in terms of the sup-norm of the forward differences also
holds.

\begin{proposition}
\label{p2.6}Let $\left( H;\left\langle \cdot ,\cdot \right\rangle \right) $
be an inner product space and $x_{i}\in H,$ $i\in \left\{ 1,\dots ,n\right\}
.$ Then we have the inequality%
\begin{equation}
\left( \sum_{i=1}^{n}\left\Vert x_{i}\right\Vert \right) ^{2}\leq \left\Vert
\sum_{i=1}^{n}x_{i}\right\Vert ^{2}+\frac{n^{2}\left( n^{2}-1\right) }{12}%
\max_{1\leq k\leq n-1}\left\Vert \Delta x_{k}\right\Vert ^{2}.  \label{2.12}
\end{equation}%
The constant $\frac{1}{12}$ is best possible in (\ref{2.12}).
\end{proposition}

\begin{proof}
As above, we have that%
\begin{equation*}
\left\Vert x_{j}-x_{i}\right\Vert \leq \sum_{k=i}^{j-1}\left\Vert \Delta
x_{k}\right\Vert \leq \left( j-i\right) \max_{1\leq k\leq n-1}\left\Vert
\Delta x_{k}\right\Vert ,
\end{equation*}%
for $1\leq i<j\leq n.$

Squaring the inequality, we get%
\begin{equation*}
\left\Vert x_{j}\right\Vert ^{2}+\left\Vert x_{i}\right\Vert ^{2}\leq 2\func{%
Re}\left\langle x_{i},x_{j}\right\rangle +\left( j-i\right) ^{2}\max_{1\leq
k\leq n-1}\left\Vert \Delta x_{k}\right\Vert ^{2}
\end{equation*}%
for any $i,j$ with $1\leq i<j\leq n,$ and since%
\begin{equation*}
2\left\Vert x_{i}\right\Vert \left\Vert x_{j}\right\Vert \leq \left\Vert
x_{j}\right\Vert ^{2}+\left\Vert x_{i}\right\Vert ^{2},
\end{equation*}%
hence%
\begin{equation}
0\leq \left\Vert x_{i}\right\Vert \left\Vert x_{j}\right\Vert -\func{Re}%
\left\langle x_{i},x_{j}\right\rangle \leq \frac{1}{2}\left( j-i\right)
^{2}\max_{1\leq k\leq n-1}\left\Vert \Delta x_{k}\right\Vert ^{2}
\label{2.13}
\end{equation}%
for any $i,j$ with $1\leq i<j\leq n.$

Applying Lemma \ref{l2.1} for $k_{ij}:=\frac{1}{2}\left( j-i\right)
^{2}\max\limits_{1\leq k\leq n-1}\left\Vert \Delta x_{k}\right\Vert ^{2},$
we can state that%
\begin{equation*}
\left( \sum_{i=1}^{n}\left\Vert x_{i}\right\Vert \right) ^{2}\leq \left\Vert
\sum_{i=1}^{n}x_{i}\right\Vert ^{2}+\sum_{1\leq i<j\leq n}\left( j-i\right)
^{2}\max_{1\leq k\leq n-1}\left\Vert \Delta x_{k}\right\Vert ^{2}.
\end{equation*}%
However,%
\begin{align*}
\sum_{1\leq i<j\leq n}\left( j-i\right) ^{2}& =\frac{1}{2}%
\sum_{i,j=1}^{n}\left( j-i\right) ^{2}=n\sum_{k=1}^{n}k^{2}-\left(
\sum_{k=1}^{n}k\right) ^{2} \\
& =\frac{n^{2}\left( n^{2}-1\right) }{12}
\end{align*}%
giving the desired inequality.

To prove the sharpness of the constant, assume that (\ref{2.12}) holds with
a constant $D>0,$ i.e.,%
\begin{equation}
\left( \sum_{i=1}^{n}\left\Vert x_{i}\right\Vert \right) ^{2}\leq \left\Vert
\sum_{i=1}^{n}x_{i}\right\Vert ^{2}+Dn^{2}\left( n^{2}-1\right) \max_{1\leq
k\leq n-1}\left\Vert \Delta x_{k}\right\Vert ^{2}  \label{2.14}
\end{equation}%
for $n\geq 2,$ $x_{i}\in H,$ $i\in \left\{ 1,\dots ,n\right\} .$

If in (\ref{2.14}) we choose $n=2,$ $x_{1}=-\frac{1}{2}e,$ $x_{2}=\frac{1}{2}%
e,$ $e\in H,$ $\left\Vert e\right\Vert =1,$ then we get $1\leq 12D$ giving $%
D\geq \frac{1}{12}.$
\end{proof}

The following result may be stated as well.

\begin{proposition}
\label{p2.7}Let $\left( H;\left\langle \cdot ,\cdot \right\rangle \right) $
be an inner product space and $x_{i}\in H,$ $i\in \left\{ 1,\dots ,n\right\}
.$ Then we have the inequality:%
\begin{equation}
\left( \sum_{i=1}^{n}\left\Vert x_{i}\right\Vert \right) ^{2}\leq \left\Vert
\sum_{i=1}^{n}x_{i}\right\Vert ^{2}+\sum_{1\leq i<j\leq n}\left( j-i\right)
^{\frac{2}{q}}\left( \sum_{k=1}^{n-1}\left\Vert \Delta x_{k}\right\Vert
^{p}\right) ^{\frac{2}{p}},  \label{2.15}
\end{equation}%
where $p>1,$ $\frac{1}{p}+\frac{1}{q}=1.$

The constant $E=1$ in front of the double sum cannot generally be replaced
by a smaller constant.
\end{proposition}

\begin{proof}
Using H\"{o}lder's inequality, we have%
\begin{align*}
\left\Vert x_{j}-x_{i}\right\Vert & \leq \sum_{k=i}^{j-1}\left\Vert \Delta
x_{k}\right\Vert \leq \left( j-i\right) ^{\frac{1}{q}}\left(
\sum_{k=i}^{j-1}\left\Vert \Delta x_{k}\right\Vert ^{p}\right) ^{\frac{1}{p}}
\\
& \leq \left( j-i\right) ^{\frac{1}{q}}\left( \sum_{k=1}^{n-1}\left\Vert
\Delta x_{k}\right\Vert ^{p}\right) ^{\frac{1}{p}},
\end{align*}%
for $1\leq i<j\leq n.$

Squaring the previous inequality, we get%
\begin{equation*}
\left\Vert x_{j}\right\Vert ^{2}+\left\Vert x_{i}\right\Vert ^{2}\leq 2\func{%
Re}\left\langle x_{i},x_{j}\right\rangle +\left( j-i\right) ^{\frac{2}{q}%
}\left( \sum_{k=1}^{n-1}\left\Vert \Delta x_{k}\right\Vert ^{p}\right) ^{%
\frac{2}{p}},
\end{equation*}%
for $1\leq i<j\leq n.$

Utilising the same argument from the proof of Proposition \ref{p2.6}, we
deduce the desired inequality (\ref{2.15}).

Now assume that (\ref{2.15}) holds with a constant $E>0,$ i.e.,%
\begin{equation*}
\left( \sum_{i=1}^{n}\left\Vert x_{i}\right\Vert \right) ^{2}\leq \left\Vert
\sum_{i=1}^{n}x_{i}\right\Vert ^{2}+E\sum_{1\leq i<j\leq n}\left( j-i\right)
^{\frac{2}{q}}\left( \sum_{k=1}^{n-1}\left\Vert \Delta x_{k}\right\Vert
^{p}\right) ^{\frac{2}{p}},
\end{equation*}%
for $n\geq 2$ and $x_{i}\in H,$ $i\in \left\{ 1,\dots ,n\right\} ,$ $p>1,$ $%
\frac{1}{p}+\frac{1}{q}=1.$

For $n=2,$ $x_{1}=-\frac{1}{2}e,$ $x_{2}=\frac{1}{2}e,$ $\left\Vert
e\right\Vert =1,$ we get $1\leq E,$ showing the fact that the inequality (%
\ref{2.15}) is sharp.
\end{proof}

The particular case $p=q=2$ is of interest.

\begin{corollary}
\label{c2.8}Let $\left( H;\left\langle \cdot ,\cdot \right\rangle \right) $
be an inner product space and $x_{i}\in H,$ $i\in \left\{ 1,\dots ,n\right\}
.$ Then we have the inequality:%
\begin{equation}
\left( \sum_{i=1}^{n}\left\Vert x_{i}\right\Vert \right) ^{2}\leq \left\Vert
\sum_{i=1}^{n}x_{i}\right\Vert ^{2}+\frac{\left( n^{2}-1\right) n}{6}%
\sum_{k=1}^{n-1}\left\Vert \Delta x_{k}\right\Vert ^{2}.  \label{2.16}
\end{equation}%
The constant $\frac{1}{6}$ is best possible in (\ref{2.16}).
\end{corollary}

\begin{proof}
For $p=q=2,$ Proposition \ref{p2.7} provides the inequality%
\begin{equation*}
\left( \sum_{i=1}^{n}\left\Vert x_{i}\right\Vert \right) ^{2}\leq \left\Vert
\sum_{i=1}^{n}x_{i}\right\Vert ^{2}+\sum_{1\leq i<j\leq n}\left( j-i\right)
\sum_{k=1}^{n-1}\left\Vert \Delta x_{k}\right\Vert ^{2},
\end{equation*}%
and since%
\begin{align*}
\sum\limits_{1\leq i<j\leq n}\left( j-i\right) & =1+\left( 1+2\right)
+\left( 1+2+3\right) +\cdots +\left( 1+2+\cdots +n-1\right)  \\
& =\sum\limits_{k=1}^{n-1}\left( 1+2+\cdots +k\right)
=\sum\limits_{k=1}^{n-1}\frac{k\left( k+1\right) }{2} \\
& =\frac{1}{2}\left[ \frac{\left( n-1\right) n\left( 2n-1\right) }{6}+\frac{%
n\left( n-1\right) }{2}\right]  \\
& =\frac{n\left( n^{2}-1\right) }{6},
\end{align*}%
hence the inequality (\ref{2.15}) is proved. The best constant may be shown
in the same way as above but we omit the details.
\end{proof}

Finally, we may state and prove the following different result.

\begin{theorem}
\label{t2.9}Let $\left( H;\left\langle \cdot ,\cdot \right\rangle \right) $
be an inner product space, $y_{i}\in H,$ $i\in \left\{ 1,\dots ,n\right\} $
and $M\geq m>0$ are such that either%
\begin{equation}
\func{Re}\left\langle My_{j}-y_{i},y_{i}-my_{j}\right\rangle \geq 0\text{ \
for }1\leq i<j\leq n,  \label{2.17}
\end{equation}%
or, equivalently,%
\begin{equation}
\left\Vert y_{i}-\frac{M+m}{2}y_{j}\right\Vert \leq \frac{1}{2}\left(
M-m\right) \left\Vert y_{j}\right\Vert \text{ \ for }1\leq i<j\leq n.
\label{2.18}
\end{equation}%
Then we have the inequality%
\begin{equation}
\left( \sum_{i=1}^{n}\left\Vert y_{i}\right\Vert \right) ^{2}\leq \left\Vert
\sum_{i=1}^{n}y_{i}\right\Vert ^{2}+\frac{1}{2}\cdot \frac{\left( M-m\right)
^{2}}{M+m}\sum\limits_{k=1}^{n-1}k\left\Vert y_{k+1}\right\Vert ^{2}.
\label{2.19}
\end{equation}%
The case of equality holds in (\ref{2.19}) if and only if%
\begin{equation}
\left\Vert y_{i}\right\Vert \left\Vert y_{j}\right\Vert -\func{Re}%
\left\langle y_{i},y_{j}\right\rangle =\frac{1}{4}\cdot \frac{\left(
M-m\right) ^{2}}{M+m}\left\Vert y_{j}\right\Vert ^{2}  \label{2.20}
\end{equation}%
for each $i,j$ with $1\leq i<j\leq n.$
\end{theorem}

\begin{proof}
Firstly, observe that, in an inner product space $\left( H;\left\langle
\cdot ,\cdot \right\rangle \right) $ and for $x,z,Z\in H,$ the following
statements are equivalent:

\begin{enumerate}
\item[(i)] $\func{Re}\left\langle Z-x,x-z\right\rangle \geq 0$

\item[(ii)] $\left\Vert x-\frac{Z+z}{2}\right\Vert \leq \frac{1}{2}%
\left\Vert Z-z\right\Vert .$
\end{enumerate}

This shows that (\ref{2.17}) and (\ref{2.18}) are obviously equivalent.

Now, taking the square in (\ref{2.18}), we get%
\begin{equation*}
\left\Vert y_{i}\right\Vert ^{2}+\frac{\left( M-m\right) ^{2}}{M+m}%
\left\Vert y_{j}\right\Vert ^{2}\leq 2\func{Re}\left\langle y_{i},\frac{M+m}{%
2}y_{j}\right\rangle +\frac{1}{n}\left( M-m\right) ^{2}\left\Vert
y_{j}\right\Vert ^{2}
\end{equation*}%
for $1\leq i<j\leq n,$ and since, obviously,%
\begin{equation*}
2\left( \frac{M+m}{2}\right) \left\Vert y_{i}\right\Vert \left\Vert
y_{j}\right\Vert \leq \left\Vert y_{i}\right\Vert ^{2}+\frac{\left(
M-m\right) ^{2}}{M+m}\left\Vert y_{j}\right\Vert ^{2},
\end{equation*}%
hence%
\begin{equation*}
2\left( \frac{M+m}{2}\right) \left\Vert y_{i}\right\Vert \left\Vert
y_{j}\right\Vert \leq 2\func{Re}\left\langle y_{i},\frac{M+m}{2}%
y_{j}\right\rangle +\frac{1}{n}\left( M-m\right) ^{2}\left\Vert
y_{j}\right\Vert ^{2},
\end{equation*}%
giving the much simpler inequality%
\begin{equation}
\left\Vert y_{i}\right\Vert \left\Vert y_{j}\right\Vert -\func{Re}%
\left\langle y_{i},y_{j}\right\rangle \leq \frac{1}{4}\cdot \frac{\left(
M-m\right) ^{2}}{M+m}\left\Vert y_{j}\right\Vert ^{2},  \label{2.21}
\end{equation}%
for $1\leq i<j\leq n.$

Applying Lemma \ref{l2.1} for $k_{ij}:=\frac{1}{4}\cdot \frac{\left(
M-m\right) ^{2}}{M+m}\left\Vert y_{j}\right\Vert ^{2},$ we deduce%
\begin{equation}
\left( \sum_{i=1}^{n}\left\Vert y_{i}\right\Vert \right) ^{2}\leq \left\Vert
\sum_{i=1}^{n}y_{i}\right\Vert ^{2}+\frac{1}{2}\frac{\left( M-m\right) ^{2}}{%
M+m}\sum_{1\leq i<j\leq n}\left\Vert y_{j}\right\Vert ^{2}  \label{2.22}
\end{equation}%
with equality if and only if (\ref{2.21}) holds for each $i,j$ with $1\leq
i<j\leq n.$

Since%
\begin{align*}
\sum_{1\leq i<j\leq n}\left\Vert y_{j}\right\Vert ^{2}& =\sum_{1<j\leq
n}\left\Vert y_{j}\right\Vert ^{2}+\sum_{2<j\leq n}\left\Vert
y_{j}\right\Vert ^{2}+\cdots +\sum_{n-1<j\leq n}\left\Vert y_{j}\right\Vert
^{2} \\
& =\sum_{j=2}^{n}\left\Vert y_{j}\right\Vert ^{2}+\sum_{j=3}^{n}\left\Vert
y_{j}\right\Vert ^{2}+\cdots +\sum_{j=n-1}^{n}\left\Vert y_{j}\right\Vert
^{2}+\left\Vert y_{n}\right\Vert ^{2} \\
& =\sum_{j=2}^{n}\left( j-1\right) \left\Vert y_{j}\right\Vert
^{2}=\sum\limits_{k=1}^{n-1}k\left\Vert y_{k+1}\right\Vert ^{2},
\end{align*}%
hence the inequality (\ref{2.19}) is obtained.
\end{proof}

\section{Further Quadratic Refinements of the Triangle Inequality}

The following lemma is of interest in itself as well.

\begin{lemma}
\label{l3.1}Let $\left( H;\left\langle \cdot ,\cdot \right\rangle \right) $
be an inner product space over the real or complex number field $\mathbb{K}$%
, $x_{i}\in H,$ $i\in \left\{ 1,\dots ,n\right\} $ and $k\geq 1$ with the
property that: 
\begin{equation}
\left\Vert x_{i}\right\Vert \left\Vert x_{j}\right\Vert \leq k\func{Re}%
\left\langle x_{i},x_{j}\right\rangle ,  \label{3.1}
\end{equation}%
for each $i,j$ with $1\leq i<j\leq n.$ Then%
\begin{equation}
\left( \sum_{i=1}^{n}\left\Vert x_{i}\right\Vert \right) ^{2}+\left(
k-1\right) \sum_{i=1}^{n}\left\Vert x_{i}\right\Vert ^{2}\leq k\left\Vert
\sum_{i=1}^{n}x_{i}\right\Vert ^{2}.  \label{3.2}
\end{equation}%
The equality holds in (\ref{3.2}) if and only if it holds in (\ref{3.1}) for
each $i,j$ with $1\leq i<j\leq n.$
\end{lemma}

\begin{proof}
Firstly, let us observe that the following identity holds true:%
\begin{align}
& \left( \sum_{i=1}^{n}\left\Vert x_{i}\right\Vert \right) ^{2}-k\left\Vert
\sum_{i=1}^{n}x_{i}\right\Vert ^{2}  \label{3.3} \\
& =\sum_{i,j=1}^{n}\left\Vert x_{i}\right\Vert \left\Vert x_{j}\right\Vert
-k\left\langle \sum_{i=1}^{n}x_{i},\sum_{j=1}^{n}x_{j}\right\rangle  \notag
\\
& =\sum_{i,j=1}^{n}\left[ \left\Vert x_{i}\right\Vert \left\Vert
x_{j}\right\Vert -k\func{Re}\left\langle x_{i},x_{j}\right\rangle \right] 
\notag \\
& =2\sum_{1\leq i<j\leq n}\left[ \left\Vert x_{i}\right\Vert \left\Vert
x_{j}\right\Vert -k\func{Re}\left\langle x_{i},x_{j}\right\rangle \right]
+\left( 1-k\right) \sum_{i=1}^{n}\left\Vert x_{i}\right\Vert ^{2},  \notag
\end{align}%
since, obviously, $\func{Re}\left\langle x_{i},x_{j}\right\rangle =\func{Re}%
\left\langle x_{j},x_{i}\right\rangle $ for any $i,j\in \left\{ 1,\dots
,n\right\} .$

Using the assumption (\ref{3.1}), we obtain%
\begin{equation*}
\sum_{1\leq i<j\leq n}\left[ \left\Vert x_{i}\right\Vert \left\Vert
x_{j}\right\Vert -k\func{Re}\left\langle x_{i},x_{j}\right\rangle \right]
\leq 0
\end{equation*}%
and thus, from (\ref{3.3}), we deduce the desired inequality (\ref{3.2}).

The case of equality is obvious by the identity (\ref{3.3}) and we omit the
details.
\end{proof}

\begin{remark}
\label{r3.4}The inequality (\ref{3.2}) provides the following reverse of the
quadratic generalised triangle inequality:%
\begin{equation}
0\leq \left( \sum_{i=1}^{n}\left\Vert x_{i}\right\Vert \right)
^{2}-\sum_{i=1}^{n}\left\Vert x_{i}\right\Vert ^{2}\leq k\left[ \left\Vert
\sum_{i=1}^{n}x_{i}\right\Vert ^{2}-\sum_{i=1}^{n}\left\Vert
x_{i}\right\Vert ^{2}\right] .  \label{3.4}
\end{equation}
\end{remark}

\begin{remark}
\label{r3.5}Since $k=1$ and $\sum_{i=1}^{n}\left\Vert x_{i}\right\Vert
^{2}\geq 0,$ hence by (\ref{3.2}) one may deduce the following reverse of
the triangle inequality%
\begin{equation}
\sum_{i=1}^{n}\left\Vert x_{i}\right\Vert \leq \sqrt{k}\left\Vert
\sum_{i=1}^{n}x_{i}\right\Vert ,  \label{3.5}
\end{equation}%
provided (\ref{3.1}) holds true for $1\leq i<j\leq n.$
\end{remark}

The following corollary providing a better bound for $\sum_{i=1}^{n}\left%
\Vert x_{i}\right\Vert ,$ holds.

\begin{corollary}
\label{c3.6}With the assumptions in Lemma \ref{l3.1}, one has the inequality:%
\begin{equation}
\sum_{i=1}^{n}\left\Vert x_{i}\right\Vert \leq \sqrt{\frac{nk}{n+k-1}}%
\left\Vert \sum_{i=1}^{n}x_{i}\right\Vert .  \label{3.6}
\end{equation}
\end{corollary}

\begin{proof}
Using the Cauchy-Bunyakovsky-Schwarz inequality%
\begin{equation*}
n\sum_{i=1}^{n}\left\Vert x_{i}\right\Vert ^{2}\geq \left(
\sum_{i=1}^{n}\left\Vert x_{i}\right\Vert \right) ^{2}
\end{equation*}%
we get%
\begin{equation}
\left( k-1\right) \sum_{i=1}^{n}\left\Vert x_{i}\right\Vert ^{2}+\left(
\sum_{i=1}^{n}\left\Vert x_{i}\right\Vert \right) ^{2}\geq \left( \frac{k-1}{%
n}+1\right) \left( \sum_{i=1}^{n}\left\Vert x_{i}\right\Vert \right) ^{2}.
\label{3.7}
\end{equation}%
Consequently, by (\ref{3.7}) and (\ref{3.2}) we deduce%
\begin{equation*}
k\left\Vert \sum_{i=1}^{n}x_{i}\right\Vert ^{2}\geq \frac{n+k-1}{n}\left(
\sum_{i=1}^{n}\left\Vert x_{i}\right\Vert \right) ^{2}
\end{equation*}%
giving the desired inequality (\ref{3.6}).
\end{proof}

The following result may be stated as well.

\begin{theorem}
\label{t3.7}Let $\left( H;\left\langle \cdot ,\cdot \right\rangle \right) $
be an inner product space and $x_{i}\in H\backslash \left\{ 0\right\} ,$ $%
i\in \left\{ 1,\dots ,n\right\} ,$ $\rho \in \left( 0,1\right) ,$ such that%
\begin{equation}
\left\Vert x_{i}-\frac{x_{j}}{\left\Vert x_{j}\right\Vert }\right\Vert \leq
\rho \ \ \text{ \ for }1\leq i<j\leq n.  \label{3.8}
\end{equation}%
Then we have the inequality%
\begin{equation}
\sqrt{1-\rho ^{2}}\left( \sum_{i=1}^{n}\left\Vert x_{i}\right\Vert \right)
^{2}+\left( 1-\sqrt{1-\rho ^{2}}\right) \sum_{i=1}^{n}\left\Vert
x_{i}\right\Vert ^{2}\leq \left\Vert \sum_{i=1}^{n}x_{i}\right\Vert ^{2}.
\label{3.9}
\end{equation}%
The case of equality holds in (\ref{3.9}) iff%
\begin{equation}
\left\Vert x_{i}\right\Vert \left\Vert x_{j}\right\Vert =\frac{1}{\sqrt{%
1-\rho ^{2}}}\func{Re}\left\langle x_{i},x_{j}\right\rangle  \label{3.10}
\end{equation}%
for any $1\leq i<j\leq n.$
\end{theorem}

\begin{proof}
The condition (\ref{3.1}) is obviously equivalent to%
\begin{equation*}
\left\Vert x_{i}\right\Vert ^{2}+1-\rho ^{2}\leq 2\func{Re}\left\langle
x_{i},\frac{x_{j}}{\left\Vert x_{j}\right\Vert }\right\rangle
\end{equation*}%
for each $1\leq i<j\leq n.$

Dividing by $\sqrt{1-\rho ^{2}}>0,$ we deduce%
\begin{equation}
\frac{\left\Vert x_{i}\right\Vert ^{2}}{\sqrt{1-\rho ^{2}}}+\sqrt{1-\rho ^{2}%
}\leq \frac{2}{\sqrt{1-\rho ^{2}}}\func{Re}\left\langle x_{i},\frac{x_{j}}{%
\left\Vert x_{j}\right\Vert }\right\rangle ,  \label{3.11}
\end{equation}%
for $1\leq i<j\leq n.$

On the other hand, by the elementary inequality%
\begin{equation}
\frac{p}{\alpha }+q\alpha \geq 2\sqrt{pq},\ \ \ p,q\geq 0,\ \alpha >0
\label{3.12}
\end{equation}%
we have%
\begin{equation}
2\left\Vert x_{i}\right\Vert \leq \frac{\left\Vert x_{i}\right\Vert ^{2}}{%
\sqrt{1-\rho ^{2}}}+\sqrt{1-\rho ^{2}}.  \label{3.13}
\end{equation}%
Making use of (\ref{3.11}) and (\ref{3.13}), we deduce that 
\begin{equation*}
\left\Vert x_{i}\right\Vert \left\Vert x_{j}\right\Vert \leq \frac{1}{\sqrt{%
1-\rho ^{2}}}\func{Re}\left\langle x_{i},x_{j}\right\rangle
\end{equation*}%
for $1\leq i<j\leq n.$

Now, applying Lemma \ref{l2.1} for $k=\frac{1}{\sqrt{1-\rho ^{2}}},$ we
deduce the desired result.
\end{proof}

\begin{remark}
\label{r3.8}If we assume that $\left\Vert x_{i}\right\Vert =1,$ $i\in
\left\{ 1,\dots ,n\right\} ,$ satisfying the simpler condition%
\begin{equation}
\left\Vert x_{j}-x_{i}\right\Vert \leq \rho \ \ \text{ \ for }1\leq i<j\leq
n,  \label{3.14}
\end{equation}%
then, from (\ref{3.9}), we deduce the following lower bound for $\left\Vert
\sum_{i=1}^{n}x_{i}\right\Vert ,$ namely%
\begin{equation}
\left[ n+n\left( n-1\right) \sqrt{1-\rho ^{2}}\right] ^{\frac{1}{2}}\leq
\left\Vert \sum_{i=1}^{n}x_{i}\right\Vert .  \label{3.15}
\end{equation}%
The equality holds in (\ref{3.15}) iff $\sqrt{1-\rho ^{2}}=\func{Re}%
\left\langle x_{i},x_{j}\right\rangle $ for $1\leq i<j\leq n.$
\end{remark}

\begin{remark}
\label{r3.9}Under the hypothesis of Proposition \ref{p2.7}, we have the
coarser but simpler reverse of the triangle inequality%
\begin{equation}
\sqrt[4]{1-\rho ^{2}}\sum_{i=1}^{n}\left\Vert x_{i}\right\Vert \leq
\left\Vert \sum_{i=1}^{n}x_{i}\right\Vert .  \label{3.16}
\end{equation}%
Also, applying Corollary \ref{c3.6} for $k=\frac{1}{\sqrt{1-\rho ^{2}}},$ we
can state that%
\begin{equation}
\sum_{i=1}^{n}\left\Vert x_{i}\right\Vert \leq \sqrt{\frac{n}{n\sqrt{1-\rho
^{2}}+1-\sqrt{1-\rho ^{2}}}}\left\Vert \sum_{i=1}^{n}x_{i}\right\Vert ,
\label{3.17}
\end{equation}%
provided $x_{i}\in H$ satisfy (\ref{3.8}) for $1\leq i<j\leq n.$
\end{remark}

In the same manner, we can state and prove the following reverse of the
quadratic generalised triangle inequality.

\begin{theorem}
\label{t3.10}Let $\left( H;\left\langle \cdot ,\cdot \right\rangle \right) $
be an inner product space over the real or complex number field $\mathbb{K}$%
, $x_{i}\in H,$ $i\in \left\{ 1,\dots ,n\right\} $ and $M\geq m>0$ such that
either%
\begin{equation}
\func{Re}\left\langle Mx_{j}-x_{i},x_{i}-mx_{j}\right\rangle \geq 0\text{ \
for }1\leq i<j\leq n,  \label{3.18}
\end{equation}%
or, equivalently, 
\begin{equation}
\left\Vert x_{i}-\frac{M+m}{2}x_{j}\right\Vert \leq \frac{1}{2}\left(
M-m\right) \left\Vert x_{j}\right\Vert \text{ \ for }1\leq i<j\leq n
\label{3.19}
\end{equation}%
hold. Then%
\begin{equation}
\frac{2\sqrt{mM}}{M+m}\left( \sum_{i=1}^{n}\left\Vert x_{i}\right\Vert
\right) ^{2}+\frac{\left( \sqrt{M}-\sqrt{m}\right) ^{2}}{M+m}%
\sum_{i=1}^{n}\left\Vert x_{i}\right\Vert ^{2}\leq \left\Vert
\sum_{i=1}^{n}x_{i}\right\Vert ^{2}.  \label{3.20}
\end{equation}%
The case of equality holds in (\ref{3.20}) if and only if%
\begin{equation}
\left\Vert x_{i}\right\Vert \left\Vert x_{j}\right\Vert =\frac{M+m}{2\sqrt{mM%
}}\func{Re}\left\langle x_{i},x_{j}\right\rangle \text{ \ for }1\leq i<j\leq
n.  \label{3.21}
\end{equation}
\end{theorem}

\begin{proof}
From (\ref{3.18}), observe that%
\begin{equation}
\left\Vert x_{i}\right\Vert ^{2}+Mm\left\Vert x_{j}\right\Vert ^{2}\leq
\left( M+m\right) \func{Re}\left\langle x_{i},x_{j}\right\rangle ,
\label{3.22}
\end{equation}%
for $1\leq i<j\leq n.$ Dividing (\ref{3.22}) by $\sqrt{mM}>0,$ we deduce%
\begin{equation*}
\frac{\left\Vert x_{i}\right\Vert ^{2}}{\sqrt{mM}}+\sqrt{mM}\left\Vert
x_{j}\right\Vert ^{2}\leq \frac{M+m}{\sqrt{mM}}\func{Re}\left\langle
x_{i},x_{j}\right\rangle ,
\end{equation*}%
and since, obviously%
\begin{equation*}
2\left\Vert x_{i}\right\Vert \left\Vert x_{j}\right\Vert \leq \frac{%
\left\Vert x_{i}\right\Vert ^{2}}{\sqrt{mM}}+\sqrt{mM}\left\Vert
x_{j}\right\Vert ^{2}
\end{equation*}%
hence%
\begin{equation*}
\left\Vert x_{i}\right\Vert \left\Vert x_{j}\right\Vert \leq \frac{M+m}{2%
\sqrt{mM}}\func{Re}\left\langle x_{i},x_{j}\right\rangle ,\text{ \ for }%
1\leq i<j\leq n.
\end{equation*}%
Applying Lemma \ref{l3.1} for $k=\frac{M+m}{2\sqrt{mM}}\geq 1,$ we deduce
the desired result.
\end{proof}

\begin{remark}
\label{r3.11}We also must note that a simpler but coarser inequality that
can be obtained from (\ref{3.20}) is 
\begin{equation*}
\left( \frac{2\sqrt{mM}}{M+m}\right) ^{\frac{1}{2}}\sum_{i=1}^{n}\left\Vert
x_{i}\right\Vert \leq \left\Vert \sum_{i=1}^{n}x_{i}\right\Vert ,
\end{equation*}%
provided (\ref{3.18}) holds true.
\end{remark}

Finally, a different result related to the generalised triangle inequality
is incorporated in the following theorem.

\begin{theorem}
\label{t3.12}Let $\left( H;\left\langle \cdot ,\cdot \right\rangle \right) $
be an inner product space over $\mathbb{K}$, $\eta >0$ and $x_{i}\in H,$ $%
i\in \left\{ 1,\dots ,n\right\} $ with the property that%
\begin{equation}
\left\Vert x_{j}-x_{i}\right\Vert \leq \eta <\left\Vert x_{j}\right\Vert 
\text{ \ for each \ }i,j\in \left\{ 1,\dots ,n\right\} .  \label{3.23}
\end{equation}%
Then we have the following reverse of the triangle inequality%
\begin{equation}
\sum_{i=1}^{n}\sqrt{\left\Vert x_{i}\right\Vert ^{2}-\eta ^{2}}\leq \frac{%
\left\Vert \sum_{i=1}^{n}x_{i}\right\Vert ^{2}}{\sum_{i=1}^{n}\left\Vert
x_{i}\right\Vert }.  \label{3.24}
\end{equation}%
The equality holds in (\ref{3.24}) iff%
\begin{equation}
\left\Vert x_{i}\right\Vert \sqrt{\left\Vert x_{j}\right\Vert ^{2}-\eta ^{2}}%
=\func{Re}\left\langle x_{i},x_{j}\right\rangle \text{ \ for each \ }i,j\in
\left\{ 1,\dots ,n\right\} .  \label{3.25}
\end{equation}
\end{theorem}

\begin{proof}
From (\ref{3.23}), we have%
\begin{equation*}
\left\Vert x_{i}\right\Vert ^{2}-2\func{Re}\left\langle
x_{i},x_{j}\right\rangle +\left\Vert x_{j}\right\Vert ^{2}\leq \eta ^{2},
\end{equation*}%
giving%
\begin{equation*}
\left\Vert x_{i}\right\Vert ^{2}+\left\Vert x_{j}\right\Vert ^{2}-\eta
^{2}\leq 2\func{Re}\left\langle x_{i},x_{j}\right\rangle ,\ \ \ i,j\in
\left\{ 1,\dots ,n\right\} .
\end{equation*}%
On the other hand, 
\begin{equation*}
2\left\Vert x_{i}\right\Vert \sqrt{\left\Vert x_{j}\right\Vert ^{2}-\eta ^{2}%
}\leq \left\Vert x_{i}\right\Vert ^{2}+\left\Vert x_{j}\right\Vert ^{2}-\eta
^{2},\ \ \ i,j\in \left\{ 1,\dots ,n\right\}
\end{equation*}%
and thus%
\begin{equation*}
\left\Vert x_{i}\right\Vert \sqrt{\left\Vert x_{j}\right\Vert ^{2}-\eta ^{2}}%
\leq \func{Re}\left\langle x_{i},x_{j}\right\rangle ,\ \ \ i,j\in \left\{
1,\dots ,n\right\} .
\end{equation*}%
Summing over $i,j\in \left\{ 1,\dots ,n\right\} ,$ we deduce the desired
inequality (\ref{3.24}).

The case of equality is also obvious from the above, and we omit the details.
\end{proof}

\section{Applications for Vector-Valued Integral Inequalities}

Let $\left( H;\left\langle \cdot ,\cdot \right\rangle \right) $ be a Hilbert
space over the real or complex number field, $\left[ a,b\right] $ a compact
interval in $\mathbb{R}$ and $\eta :\left[ a,b\right] \rightarrow \lbrack
0,\infty )$ a Lebesgue integrable function on $\left[ a,b\right] $ with the
property that $\int_{a}^{b}\eta \left( t\right) dt=1.$ If, by $L_{\eta
}\left( \left[ a,b\right] ;H\right) $ we denote the Hilbert space of all
Bochner measurable functions $f:\left[ a,b\right] \rightarrow H$ with the
property that $\int_{a}^{b}\eta \left( t\right) \left\Vert f\left( t\right)
\right\Vert ^{2}dt<\infty ,$ then the norm $\left\Vert \cdot \right\Vert
_{\eta }$ of this space is generated by the inner product $\left\langle
\cdot ,\cdot \right\rangle _{\eta }:H\times H\rightarrow \mathbb{K}$ defined
by 
\begin{equation*}
\left\langle f,g\right\rangle _{\eta }:=\int_{a}^{b}\eta \left( t\right)
\left\langle f\left( t\right) ,g\left( t\right) \right\rangle dt.
\end{equation*}%
The following proposition providing a reverse of the integral generalised
triangle inequality may be stated.

\begin{proposition}
\label{p4.1} Let $\left( H;\left\langle \cdot ,\cdot \right\rangle \right) $
be a Hilbert space and $\eta :\left[ a,b\right] \rightarrow \lbrack 0,\infty
)$ as above. If $g\in L_{\eta }\left( \left[ a,b\right] ;H\right) $ is so
that $\int_{a}^{b}\eta \left( t\right) \left\Vert g\left( t\right)
\right\Vert ^{2}dt=1$ and $f_{i}\in L_{\eta }\left( \left[ a,b\right]
;H\right) ,i\in \left\{ 1,\dots ,n\right\} ,$ and $M\geq m>0$ are so that
either%
\begin{equation}
\func{Re}\left\langle Mf_{j}\left( t\right) -f_{i}\left( t\right)
,f_{i}\left( t\right) -mf_{j}\left( t\right) \right\rangle \geq 0
\label{4.1}
\end{equation}%
or, equivalently, 
\begin{equation*}
\left\Vert f_{i}\left( t\right) -\frac{m+M}{2}f_{j}\left( t\right)
\right\Vert \leq \frac{1}{2}\left( M-m\right) \left\Vert f_{j}\left(
t\right) \right\Vert
\end{equation*}%
for a.e. $t\in \left[ a,b\right] $ and $1\leq i<j\leq n,$ then we have the
inequality%
\begin{eqnarray}
&&\left[ \sum_{i=1}^{n}\left( \int_{a}^{b}\eta \left( t\right) \left\Vert
f_{i}\left( t\right) \right\Vert ^{2}dt\right) ^{1/2}\right] ^{2}
\label{4.2} \\
&\leq &\int_{a}^{b}\eta \left( t\right) \left\Vert \sum_{i=1}^{n}f_{i}\left(
t\right) \right\Vert ^{2}dt  \notag \\
&&+\frac{1}{2}\cdot \frac{\left( M-m\right) ^{2}}{m+M}\int_{a}^{b}\eta
\left( t\right) \left( \sum_{k=1}^{n-1}k\left\Vert f_{k+1}\left( t\right)
\right\Vert ^{2}\right) dt.  \notag
\end{eqnarray}%
The case of equality holds in (\ref{4.2}) if and only if%
\begin{eqnarray*}
&&\left( \int_{a}^{b}\eta \left( t\right) \left\Vert f_{i}\left( t\right)
\right\Vert ^{2}dt\right) ^{1/2}\left( \int_{a}^{b}\eta \left( t\right)
\left\Vert f_{j}\left( t\right) \right\Vert ^{2}dt\right) ^{1/2} \\
&&-\int_{a}^{b}\eta \left( t\right) \func{Re}\left\langle f_{i}\left(
t\right) ,f_{j}\left( t\right) \right\rangle dt \\
&=&\frac{1}{4}\cdot \frac{\left( M-m\right) ^{2}}{m+M}\int_{a}^{b}\eta
\left( t\right) \left\Vert f_{j}\left( t\right) \right\Vert ^{2}dt
\end{eqnarray*}%
for each $i,j$ with $1\leq i<j\leq n.$
\end{proposition}

\begin{proof}
We observe that 
\begin{eqnarray*}
&&\func{Re}\left\langle Mf_{j}-f_{i},f_{i}-mf_{j}\right\rangle _{\eta } \\
&=&\int_{a}^{b}\eta \left( t\right) \func{Re}\left\langle Mf_{j}\left(
t\right) -f_{i}\left( t\right) ,f_{i}\left( t\right) -mf_{j}\left( t\right)
\right\rangle dt\geq 0
\end{eqnarray*}%
for any $i,j$ with $1\leq i<j\leq n.$

Applying Theorem \ref{t2.9} for the Hilbert space $L_{\eta }\left( \left[ a,b%
\right] ;H\right) $ and for $y_{i}=f_{i},i\in \left\{ 1,\dots ,n\right\} ,$
we deduce the desired result.
\end{proof}

Another integral inequality incorporated in the following proposition holds:

\begin{proposition}
\label{p4.2} With the assumptions of Proposition \ref{p4.1}, we have 
\begin{eqnarray}
&&\frac{2\sqrt{mM}}{m+M}\left[ \sum_{i=1}^{n}\left( \int_{a}^{b}\eta \left(
t\right) \left\Vert f_{i}\left( t\right) \right\Vert ^{2}dt\right) ^{1/2}%
\right] ^{2}  \label{4.3} \\
&&+\frac{\left( \sqrt{M}-\sqrt{m}\right) ^{2}}{m+M}\sum_{i=1}^{n}%
\int_{a}^{b}\eta \left( t\right) \left\Vert f_{i}\left( t\right) \right\Vert
^{2}dt  \notag \\
&\leq &\int_{a}^{b}\eta \left( t\right) \left\Vert \sum_{i=1}^{n}f_{i}\left(
t\right) \right\Vert ^{2}dt.  \notag
\end{eqnarray}%
The case of equality holds in (\ref{4.3}) if and only if 
\begin{eqnarray*}
&&\left( \int_{a}^{b}\eta \left( t\right) \left\Vert f_{i}\left( t\right)
\right\Vert ^{2}dt\right) ^{1/2}\left( \int_{a}^{b}\eta \left( t\right)
\left\Vert f_{j}\left( t\right) \right\Vert ^{2}dt\right) ^{1/2} \\
&=&\frac{M+m}{2\sqrt{mM}}\int_{a}^{b}\eta \left( t\right) \func{Re}%
\left\langle f_{i}\left( t\right) ,f_{j}\left( t\right) \right\rangle dt
\end{eqnarray*}%
for any $i,j$ with $1\leq i<j\leq n.$
\end{proposition}

The proof is obvious by Theorem \ref{t3.10} and we omit the details.

\end{document}